\documentclass[12pt,english]{article}
\usepackage{amsmath}

\newcommand{\Cg}{\mathrm{Cg}}

\newcommand{\variedad}[1]{\mathcal{#1}}

\newcommand{\R}{\modelo{R}}

\newcommand{\V}{\variedad{V}}

\newcommand{\id}{\approx}

\renewcommand{\id}{\approx}

\renewcommand{\o}{\vee}

\newcommand{\y}{\wedge}

\usepackage{amsmath}
\usepackage{amsthm}
\usepackage{verbatim}
\usepackage{babel}
\usepackage[utf8]{inputenc}
\newlength{\ancho}
\textwidth=\ancho
\newcommand{\mb}[1]{\mathbf{#1}}

\renewcommand{\>}{\rangle}
\newcommand{\por}{\cdot}
\newcommand{\men}{\preceq}
\newcommand{\may}{\succeq}
\newcommand{\impl}{\rightarrow}
\renewcommand{\R}{\variedad{R}}
\theoremstyle{remark}

\theoremstyle{plain}
\newtheorem{theorem}{Theorem}
\newtheorem{lemma}[theorem]{Lemma}
\newtheorem{claim}[theorem]{Claim}

\theoremstyle{definition}
\newtheorem{definition}[theorem]{Definition}

\begin{document}
\title{Directly Indecomposables in Semidegenerate Varieties of Connected po-Groupoids}  
\author{Pedro S\'anchez Terraf\thanks{Supported by Conicet. \newline
}
}
\date{}
\maketitle  

\begin{abstract}
We study varieties with a term-definable poset structure, \emph{po-groupoids}. It is known
that connected posets have the \emph{strict refinement property}
(SRP). In \cite{DFC} it is proved that 
semidegenerate varieties with the SRP have definable factor congruences
and if the similarity type is finite, directly indecomposables are
axiomatizable by a set of first-order  sentences. We obtain such a set
for semidegenerate varieties of connected po-groupoids 
and show its quantifier complexity is bounded in general.
\end{abstract}

% 3515638176 Maxi
% 30

\section{Introduction and Basic Definitions}
\begin{definition}
A \emph{po-groupoid} is a groupoid $\<A,\por\>$ such that the relation
defined by 
\[x \men y  \text{ if and only if } x \por y = x\]
is a partial order on $A$, the \emph{order related to $\<A,\por\>$}. A \emph{variety of po-groupoids} is an
equational class $\V$ of algebras with binary term $\por$
in the language of $\V$ such that for every $A\in\V$, $\<A,\por^A\>$
is a po-groupoid.
\end{definition}
For every poset $\<A,\men\>$ one can define a po-groupoid operation
$*$ on $A$ setting
\[x * y  := \begin{cases} x & \text{if } x\men y \\ 
                y  & \text{if } x\not\men y. \end{cases} \]
such that $\men$ is the order related to
$\<A,*\>$.

Po-groupoids are obviously idempotent, but need not be associative
nor commutative.  Examples of po-groupoids are semilattices and, more
generally, the variety axiomatized by the following identities:
\begin{equation}\label{eq:2}
\begin{split}
(x \por y) \por z &\id x \por (y \por z) \\
x \por x &\id x \\
x \por y \por x &\id y \por x.
\end{split}
\end{equation}
This variety is exactly the class of associative po-groupoids,
\emph{po-semigroups} for short (see Claim~\ref{cl:po-semigroups}). J. Gerhard \cite{Ger2} proved it is
not residually small.  
 
A poset $\<A,\men\>$ is said to be \emph{connected} if the associated graph
is. Equivalently, $\<A,\men\>$ is connected if  for all $x,y\in A$,
there exists a positive integer $n$ and elements $m_1,\dots,m_{2n-1}$ such
that
\begin{equation}\label{eq:1}
x \may m_1 \men m_2 \may \dots \men m_{2n-2} \may m_{2n-1} \men y
\end{equation}
where $\may$ is the converse relation to $\men$. A po-groupoid will
be called connected if the related order is.

A variety $\V$ is \emph{semidegenerate} if no non trivial member has a
trivial subalgebra. Equivalently, if every $A\in\V$ has a compact
universal congruence (Kollar \cite{5}). In this work we will consider semidegenerate varieties of
connected po-groupoids over a finite language.

 It is known \cite[Section 5.6]{4}
that connected posets have  the \emph{strict refinement property}
(SRP). In  a joint work with D. Vaggione \cite{DFC}  it is proved that 
semidegenerate varieties with the SRP have definable factor
congruences and if the similarity type is finite, directly indecomposables are
axiomatizable by a set of first-order  sentences. 
Our main result is an application of \cite{DFC}
and a result of R.~Willard~\cite{7}.
\begin{theorem}\label{th:main}
Let $\V$ be a semidegenerate variety of connected po-groupoids over
a finite language. Then the 
class of directly indecomposable algebras of $\V$ is axiomatizable by
a $\Pi_6$ (i.e., $\forall\exists\forall\exists\forall\exists$)
first-order sentence plus axioms for $\V$.
\end{theorem}

%% \section{Basic Properties of Po-Groupoids}
%% We now list some basic results about po-groupoids. The proofs are
%% easy exercises and therefore they are omited.
%% \begin{prop}
%% Let  $\<A,\por\>$ be a po-groupoid and $\men$ its related order.
%% \begin{enumerate}
%% \item  $x \por y \men y$.
%% \item  Right multiplication is monotonous:  $x \men y \impl x \por z \men y \por z$.
%% \item  Comparable elements commute  $x \men y \impl x \por y = y \por x$.
%% \item  If $x$ and $y$ commute, $x\por y$ is the greatest lower bound
%%   of $\{x,y\}$ (but not conversely).
%% \end{enumerate}
%% \end{prop}

%% \begin{definition}
%% A \emph{relative semilattice} is a po-groupoid that satisfies
%% $x\por y \por z \id y \por x \por z$. 
%% \end{definition}
%% \begin{prop}
%% Relative semilattices are exactly the po-groupoids that satisfy
%% $x\men z \y y\men z \impl x\por y = y\por x$, i.e., every bounded
%% subset is a semilattice.
%% \end{prop}
%% \begin{prop}
%% Let  $\<A,\por\>$ be a relative semilattice and $\men$ its related order.
%% \begin{enumerate}
%% \item  Multiplication is monotonous.
%% \item  $x \por y = y \por x \y z \men y \impl x \por z = z \por x$.
%% \end{enumerate}
%% Each of the previous properties fail in some po-groupoid.
%% \end{prop}

\section{Results}
We first state two auxiliary results that provide a useful Mal'cev
condition for semidegeneracy. If $\vec{a},\vec{b} \in A^n$,  $\Cg
^{A}(\vec{a},\vec{b})$ will stand 
for the $A$-congruence generated by $(a_1,b_1), \dots, (a_n,b_n)$.
\begin{lemma}[Mal'cev]\label{malsev}Let $A$ be any algebra and let $a,b\in A,$
 $\vec{a},\vec{b}\in
A^{n}.$ Then $(a,b)\in \Cg ^{A}(\vec{a},\vec{b})$ if and only if 
there
exist $(n+m)$-ary terms 
$p_{1}(\vec{x},\vec{u}),\dots,p_{k}(\vec{x},\vec{u})$,
with $k$ odd and, $\vec{u}\in A^{m}$ such that: 
\begin{equation*}
\begin{split}
a& =p_{1}(\vec{a},\vec{u}) \\
p_{i}(\vec{b},\vec{u})& =p_{i+1}(\vec{b},\vec{u}),\ i\text{
odd} \\
p_{i}(\vec{a},\vec{u})& =p_{i+1}(\vec{a},\vec{u}),\ i\text{
even} \\
p_{k}(\vec{b},\vec{u})& =b
\end{split}
\end{equation*}
\end{lemma}
\begin{lemma}\label{l:0y1}
If $\V$ is semidegenerate, there exist  positive integers $l$ and  $k$
(with $k$ odd), unary terms
$0_{1}(w),\dots,0_{l}(w),$ $1_{1}(w),\dots,1_{l}(w)$ and  $(2+l)$-ary
terms $U_{i}(x,y,\vec{z}),$ $i=1,\dots,k,$ such that the following
identities hold in $\mathcal{V}$:
\begin{equation}\label{eq:34}
  \begin{split}
    x&\id U_{1}(x,y,\vec 0)\\
    U_{i}(x,y,\vec{1})&\id U_{i+1}(x,y,\vec 1) \text{ with $i$ odd}\\
    U_{i}(x,y,\vec 0)&\id U_{i+1}(x,y,\vec 0) \text{ with $i$ even}\\
    U_{k}(x,y,\vec{1})&\id y
    \end{split}
\end{equation}
where  $w$, $x$ and $y$ are distinct variables, $\vec{0}=(0_{1}(w),\dots,0_{l}(w))$ and\ $\vec{1}=(1_{1}(w),\dots,1_{l}(w)).$
\end{lemma}
\begin{proof}
By Kollar \cite{5} every algebra in a semidegenerate variety has a
compact universal congruence. By Lemma 3 in Vaggione \cite{va5}, there
are terms $0_{1}(w),\dots,0_{l}(w),$ $1_{1}(w),\dots,1_{l}(w)$ such that 
$\Cg^A(\vec 0 , \vec 1)$ is the
universal congruence in each $A\in \V$.  Now apply Lemma~\ref{malsev}.
\end{proof}

 Throughout this paper we will assume  that we may find
 closed terms $\vec{0}$
and  $\vec{1}$ for $\V$. Of course, this can be
achieved when the language has a constant symbol and we will make this
assumption in order to   clarify our treatment. The
proofs remain valid in the general case.

The next lemma ensures that we  have a uniform way to witness connection.
\begin{lemma}\label{l:conexion}
For every variety $\V$ of connected po-groupoids, there exist a
positive integer $n$ and binary terms
$m_i(x,y)$, $i=1,\dots,2n-1$ in the language of $\V$ such that  the following
identities hold in $\mathcal{V}$:
\begin{align*}
m_1(x,y) \por x &\id m_1(x,y) \\
m_i(x,y) \por m_{i\pm 1}(x,y) & \id m_i(x,y) & \text{if $i<2n-2$ is
  odd} \\
m_{2n-1}(x,y) \por y &\id m_{2n-1}(x,y).
\end{align*}
\end{lemma}
\begin{proof}
It is enough to consider the $\V$-free algebra freely generated 
by $\{x,y\}$. The elements $m_i$ that connect $x$ and $y$ are binary
terms, and the desired equations are equivalent to the assertions in~(\ref{eq:1}).
\end{proof}

For the rest of this section, $\V$ will denote a semidegenerate variety
of po-groupoids with $k$ terms $U_i$ and $2n-1$ terms $m_i$ as in
the previous lemmas. In the following we define formulas $\psi$, $\phi$ and $\pi$
 codifying the fact that connected po-groupoids have the SRP. They appear
under the same names in Willard
\cite[Section 5]{7}, though
the key improvement is a radical simplification of the first
one.

A formula is called \emph{$\V$-factorable} (\emph{factorable} for short) if it belongs to the smallest set $\mb{F}$
containing every atomic formula that is closed under conjunction,
existential and universal quantification, and the following rule:
\begin{quote}
if $\alpha(\vec x),\beta(\vec x, \vec y), \gamma(\vec x, \vec y)\in
\mb{F}$ and $\V\models\forall \vec x \bigl(\alpha(\vec x) \impl \exists \vec y
  \beta(\vec x, \vec y)\bigr)$ then $\forall \vec y   \bigl( \beta(\vec x,
  \vec y) \impl  \gamma(\vec x, \vec y)\bigr) \in \mb{F}$.
\end{quote}
Factorable formulas are preserved by direct factors and direct
products (see \cite[Section 1]{7}).
\begin{lemma}\label{l:psi-n}
There exists a factorable $\Pi_2$ formula $\psi(x,y,z)$ such that:
\begin{enumerate}
\item \label{item:1}$\V\models \psi(x,y,x)$.
\item \label{item:2}$\V\models \psi(x,y,y)$.
\item \label{item:3}$\V\models \psi(x,x,z) \impl x\men z$.
\end{enumerate}
\end{lemma}
\begin{proof}
Let $\psi(x,y,z)$ be the following formula:
\begin{multline*}
\forall u_1,\dots,u_{2n-1} \\
  u_1 \por x = u_1 \por u_2 \y \ 
\bigwedge_{i=2}^{n-1} u_{2i-1}\por u_{2i-2} =  u_{2i-1}\por u_{2i}\ \  \y
\ u_{2n-1} \por u_{2n-2} = u_{2n-1} \por y \impl \\
\impl \exists v_1,\dots,v_{n-1} : 
u_1 \por x = u_1 \por v_1 \y \ 
\bigwedge_{i=2}^{n-1} u_{2i-1}\por v_{i-1} =  u_{2i-1}\por v_{i} \  \ \y
\ 
u_{2n-1} \por v_{n-1} = u_{2n-1} \por z.
\end{multline*}
We may see that $\psi$ is factorable observing that $u_i:=m_i(x,y)$
satisfies the antecedent for any choice of $x$ and $y$.

To prove~\ref{item:1}, simply take $v_i:=x$ for all
$i$. For~\ref{item:2}, it suffices to assign $v_i := u_{2i}$ for all
$i$. Finally, suppose $\psi(x,x,z)$ holds. Take $u_i
:= x$. With this choice, the antecedent holds. The consequent turns
to:
\[x \por x = x \por v_1 \ \y \ 
\bigwedge_{i=2}^{n-1} x\por v_{i-1} =  x\por v_{i} \  \y \   
x \por v_{n-1} = x \por z,\]
from which we conclude $x = x\por z$, and we have proved~\ref{item:3}.
\end{proof}

\begin{lemma}\label{l:pi}
There exists a factorable $\Pi_3$ formula $\pi(x,y,z,w)$ such that:
\begin{enumerate}
\item \label{item:4}$\V\models \pi(x,x,z,w)$
\item \label{item:5} $\V\models \pi(x,y,x,y)$
\item \label{item:6}$\V\models \pi(x,y,z,z) \impl x=y$
\end{enumerate}
\end{lemma}
\begin{proof}
We first define $\phi(x,y,w_1,w_2)$ to be $\psi(x,w_1,y) \y
\psi(y,w_2,x)$. Take $\pi(x,y,z,w)$ to be $\forall w_1,w_2 :\phi(z,w,w_1,w_2) \impl
\phi(x,y,w_1,w_2)$.

It is immediate that \ref{item:5} holds. Property \ref{item:4} holds
thanks to Lemma~\ref{l:psi-n}(\ref{item:1}). And we can check
\ref{item:6} by taking $w_1 := y$ and $w_2 :=x$. 

Finally, $\pi$ is factorable since $\exists w_1,w_2
:\phi(z,w,w_1,w_2)$ holds in $\V$ due to
Lemma~\ref{l:psi-n}(\ref{item:2}).
\end{proof}
The formula $\Phi$ appearing in the next lemma  is a first-order definition of
factor congruences in $\V$ using central elements. This concept (in
its full generality) is due to Vaggione \cite{va0}.

If $\vec a \in A^l$ and $\vec b \in B^l$, we will write $[\vec a, \vec b]$ in place
of $ ((a_{1},b_{1}),\dots,(a_{l},b_{l})) \in (A \times B)^l$.
If $A\in \mathcal{V}$, we say that $\vec{e}\in A^{l}$ is a \emph{central element} of $A$
if there exists an isomorphism $A\rightarrow A_{1}\times A_{2}$ such
that 
\[\vec{e}\mapsto [ \vec{0},\vec{1}].\]
The set of all central elements of $A$ will be called the
\emph{center} of $A$.
Two central elements $\vec{e},\vec{f}$ will be called 
\emph{complementary} if there exists an isomorphism $ A\rightarrow A_{1}\times A_{2}$
such that $\vec{e}\mapsto [\vec{0},\vec{1}]$ and $\vec{f}%
\mapsto [\vec{1},\vec{0}]$. It is immediate that an algebra is
directly indecomposable if and only if it has exactly two  central
elements, namely 
$\vec 0$ and $\vec 1$.
\begin{lemma}
There exists a factorable $\Sigma_4$ formula $\Phi(x,y,\vec z)$ such
that for all $A, B\in \V$, and $a,c\in A$,
  $b,d\in B$,
  \[A\times B \models \Phi\bigl(\<a,b\>, \<c,d\>, [\vec 0, \vec 1]\bigr)
  \quad \text{ if and only if } \quad a=c.\]
\end{lemma}
\begin{proof}
Take  $\Phi(x,y,\vec z)$ to be
\begin{multline*}
\exists a_1,\dots,a_{n-1} :  \pi(x,a_1,U_1(x,y,\vec z),U_1(x,y,\vec w)) 
\y \bigwedge_{i\text{ odd}} \pi(a_i,a_{i+1},U_{i+1}(x,y,\vec
w),U_{i+1}(x,y,\vec z)) \y \\
\y \bigwedge_{i\text{ even}} \pi(a_i,a_{i+1},U_{i+1}(x,y,\vec
z),U_{i+1}(x,y,\vec w))
\ \ \y \   \pi(a_{n-1},y,U_k(x,y,\vec z),U_k(x,y,\vec w)).
\end{multline*}

Suppose that $A\times B \models \Phi\bigl(\<a,b\>, \<c,d\>, [\vec 0,
  \vec 1]\bigr)$. Looking at the first coordinate, we obtain
\begin{multline*}
\pi(a,a^1_1,U_1(a,c,\vec 0),U_1(a,c,\vec 0)) 
\y \bigwedge_{i\text{ odd}} \pi(a^1_i,a^1_{i+1},U_{i+1}(a,c,\vec
0),U_{i+1}(a,c,\vec 0)) \;\y \\
\y \bigwedge_{i\text{ even}} \pi(a^1_i,a^1_{i+1},U_{i+1}(a,c,\vec
0),U_{i+1}(a,c,\vec 0))
\ \y \   \pi(a^1_{n-1},c,U_k(a,c,\vec 0),U_k(a,c,\vec 0)),
\end{multline*}
where $a_i=\<a^1_i,a^2_i\>$. Applying Lemma~\ref{l:pi}(\ref{item:6})
we obtain $a=c$.

Finally, to show that 
\[A\times B \models \Phi\bigl(\<a,b\>, \<a,d\>, [\vec 0, \vec
  1]\bigr)\]
it suffices to take $a_i :=
\<a,U_i(b,d,\mathbf{i})\>$, where $\mathbf{i}$ is $\vec 0$ if
$i$ is even, otherwise $\vec 1$. This is routine.
\end{proof}
Now we define the center of algebras in $\V$ in first-order logic.
\begin{lemma}
There is a $\Pi_5$ formula $\zeta(\vec z, \vec w)$ such that  for all
$A\in\V$ and  $\vec e,\vec f \in A^l$ we have
that $\vec e$ and $\vec f$ are complementary 
central elements if and only if $A\models
\zeta(\vec e,\vec f)$. 
\end{lemma}
\begin{proof}
The following formulas in the language of $\V$ will assert the properties
needed to force $\Phi(\cdot,\cdot,\vec z)$ and  $\Phi(\cdot,\cdot,\vec w)$ to
define the pair of complementary factor congruences associated with
$\vec z$ and $\vec w$. The names are almost self-explanatory.
\begin{itemize}
\item $\mathit{CAN}(\vec z,\vec w\,)= \bigwedge_{i=1}^l \Phi(0_i,z_i,\vec z\,) \y \bigwedge_{i=1}^l \Phi(1_i,w_i,\vec z\,)$
\item $\mathit{PROD}(\vec z,\vec w\,)=\forall x,y\exists z\ 
\Bigl(\Phi(x,z,\vec z\,)\wedge \Phi(z,y,\vec w\,)\Bigr)$
\item $\mathit{INT}(\vec z,\vec w\,)=\forall x,y\ 
\Bigl(\Phi(x,y,\vec z\,)\wedge \Phi(x,y,\vec w\,)\rightarrow x=y\Bigr)$
\item $\mathit{REF}(\vec z,\vec w\,)=\forall x\ \Phi(x,x,\vec z\,)$
\item $\mathit{SYM}(\vec z,\vec w\,)=\forall x,y,z\ \Bigl(\Phi(x,y, \vec z\,)\wedge \Phi(y,z,\vec z\,)\wedge 
\Phi(z,x,\vec w\,) \rightarrow z=x\Bigr)$
\item $\mathit{TRANS}(\vec z,\vec w\,)=\forall x,y,z,u
\;\Bigl(\Phi(x,y,\vec z\,)\wedge 
\Phi(y,z,\vec z\,)\wedge \Phi(x,u,\vec z\,)\wedge \Phi(u,z,\vec w\,)\rightarrow u=z\Bigr)$
\item For each  $m$-ary function symbol $F$:
\begin{multline*}
\mathit{PRES}_{F}(\vec z,\vec w\,)=\forall u_{1},v_{1},\dots,u_{m},v_{m} \\ 
\Bigl(\bigwedge_{j}\Phi(u_{j},v_{j},\vec z\,)\Bigr)\;\wedge 
\Phi(F(u_{1},\dots,u_{m}),z,\vec z\,)\wedge 
\Phi(z,F(v_{1},\dots,v_{m}),\vec w\,)\rightarrow\\ 
\rightarrow z=F(v_{1},\dots,v_{m})
\end{multline*}
\end{itemize}
Now take  $\mathit{CAN'}$, $\mathit{REF}',$ $\mathit{SYM}',$ $\mathit{TRANS}'$ and $\mathit{PRES}_{F}'$ to 
be
the result of interchanging $\vec z$ with $\vec w$ in $CAN,$ $REF,$ $SYM,$ 
$TRANS$
and $PRES_{F},$ respectively, and let $\zeta$ be the
  conjunction of:
\begin{gather*}
\bigwedge\{\mathit{CAN},\mathit{PROD},\mathit{INT},
\mathit{REF},\mathit{SYM},\mathit{TRANS},\mathit{CAN}',\mathit{REF}',\mathit{SYM}',\mathit{TRANS}'\}\\
\bigwedge\{PRES_{F},\mathit{PRES}_{F}' : F\text{ a function
  symbol}\}.
\end{gather*}
Details can be found in \cite[Lemma 4.1]{DFC}.
\end{proof}

\begin{proof}[Proof of Theorem~\ref{th:main}]
The formula
\[\vec 0 \neq \vec 1 \ \y \ \forall \vec e, \vec f : \zeta(\vec e,\vec f) \ \impl\; \bigl((\vec e =
\vec 0 \y \vec f= \vec 1) \o
(\vec e = \vec 1 \y \vec f=\vec 0)\bigr)\]
together with axioms for $\V$ defines the subclass of directly indecomposables.
\end{proof}

\section{Examples}
We first mention that we cannot eliminate the semidegeneracy
hypothesis, since even the class of directly indecomposable lattices
with 0 is not axiomatizable in first-order logic (see Willard \cite{indec-lat}). Also, in
\cite[Section 6]{DFC} it is shown that semidegeneracy by itself does
not ensure definability of directly indecomposables. By considering in
this last case a trivial (antichain) po-groupoid structure (for instance, defining
$x\por y := y$ for all $x,y$)  we
deduce that an arbitrary semidegenerate variety of po-groupoids may
not have a first-order-axiomatizable class of indecomposables; hence
we cannot drop connectedness.

We will now consider the following variety $\R$:
%%  (hereafter we omit the
%% ``$\por$''):
%% \renewcommand{\por}{\,}
\begin{equation}\label{eq:rel_semilatt}
\begin{split}
(x \por y) \por z &\id x \por (y \por z) \\
x \por x &\id x \\
x \por y \por z &\id y \por x \por z.
\end{split}
\end{equation}
This is obviously a variety of po-groupoids. Moreover, the variety of
po-semigroups (defined by equations~(\ref{eq:2})) covers $\R$ in the
lattice  of  equational classes of idempotent semigroups (see 
Gerhard \cite{Ger1}).

From the partial-order point of view, this groupoids are ``relative
meet-semilattices'': whenever $A\in\V$, every bounded subalgebra of
$A$ is a semilattice. Actually, the third axiom is equivalent to this
property.
\begin{lemma}\label{l:rel_semilatt}
Let $\V$ be a variety of associative po-groupoids. The following are
equivalent:
\begin{enumerate}
\item $\V\models x\men z \y y\men z \impl x\por y = y\por x$.
\item $\V\models x \por y \por z \id y \por x \por z$.
\end{enumerate}
\end{lemma}
\begin{proof}
We first need
\begin{claim}\label{cl:po-semigroups}
Every associative po-groupoid satisfies $x \por y \por x \id y \por x$.
\end{claim}
\begin{proof}
We now that $x\por y \men y$. We obtain immediately that  $x \por y
\por x \men y \por x$ and $y \por x = y \por x \por y \por x \men   x
\por y \por x$. By antisymmetry we get $x \por y \por x = y \por x$.
\end{proof}
Assume (1). Note that $x\por y \por z\men z$ and $y \por x \por z \men
z$, hence $x\por y \por z \por y \por x \por z = y \por x \por z \por
x\por y \por z$. We may simplify this expression using the Claim to obtain (2).

Now suppose (2) holds, and assume $ x, y\men z$. Hence
\[x \por y = x \por y \por z = y \por x \por z = y \por x,\]
and we have (1).
\end{proof}

\begin{lemma}\label{l:cota_inf}
Let $\V$ will be a variety of connected
po-groupoids that satisfies the defining
identities of $\R$, and let $x_1,\dots,x_j \in A\in\V$. Then $A\models \exists u :
\bigwedge_{i=2}^j u \por x_1
= u \por x_i$. 
\end{lemma}
\begin{proof}
We only prove the case $j=2$, from which the rest
can be easily derived.
We will prove by induction on $n$:
\begin{quote}
For all  $m_1,\dots,m_{2n-1} \in A$ such
that
\begin{equation*}
 m_1 \men m_2 \may \dots \men m_{2n-2} \may m_{2n-1} 
\end{equation*}
we have $m_1 \por m_3 \por \dots \por  m_{2n-1} =  m_{2n-1} \por \dots
\por m_3 \por m_1$.
\end{quote}
If $n=2$, we have $m_1, m_3 \men m_2$, hence by
Lemma~\ref{l:rel_semilatt}  $m_1 \por m_3 =  m_3 \por m_1$. Now
suppose the assertion holds for $n$, and assume $m_1,\dots,m_{2n+1}
\in A$ satisfy:
\begin{equation*}
 m_1 \men m_2 \may \dots \men m_{2n-2} \may m_{2n-1} \men m_{2n} \may m_{2n+1} 
\end{equation*}
We may apply the inductive hypothesis and obtain
\begin{equation*}
m_1 \por m_3 \por \dots \por  m_{2n-1}\por m_{2n+1} = m_1 \por
m_{2n+1}\por  m_{2n-1} \por \dots \por m_3 
\end{equation*}
This last term equals $m_{2n+1}\por  m_{2n-1} \por \dots \por
m_1 \por m_3$ (by the third axiom of $\R$) and by the case $n=2$ we
may commute $m_1$ and $m_3$, obtaining:
\begin{equation}\label{eq:3}
m_1 \por m_3 \por \dots \por  m_{2n-1}\por m_{2n+1} = 
m_{2n+1}\por  m_{2n-1} \por \dots \por m_3 \por m_1.
\end{equation}
Once we have this, we may take $u := m_1 \por m_3 \por \dots \por  m_{2n-1}$, where
\begin{equation*}
x_1 \may  m_1 \men m_2 \may \dots \men m_{2n-2} \may m_{2n-1} \men x_2 
\end{equation*}
since $u \por x_2 = m_1 \por m_3 \por \dots \por  m_{2n-1} \por x_2 =  m_1 \por m_3
\por \dots \por  m_{2n-1} = m_{2n-1} \por \dots
\por m_3 \por m_1 = m_{2n-1} \por \dots
\por m_3 \por m_1 \por x_1 = u \por x_1$
\end{proof}
The previous lemmas were discovered using
the \emph{Prover9/Mace4} program bundle by W.~McCune
\cite{p9m4,Mace4}.

\begin{lemma}
Let $\V$ will be a semidegenerate variety of connected
po-groupoids that satisfies the defining
identities of $\R$. There exists a factorable $\Pi_1$ formula $\Phi(x,y,\vec z)$ such
that for all $A, B\in \V$, and $a,c\in A$,
  $b,d\in B$,
  \[A\times B \models \Phi\bigl(\<a,b\>, \<c,d\>, [\vec 0, \vec 1]\bigr)
  \quad \text{ if and only if } \quad a=c.\]
\end{lemma}
\begin{proof}
To fix notation, we assume again there are $k$ terms
$U_i$ and $2n-1$ terms $m_i$ that witness semidegeneracy and
connection, respectively (see Lemmas~\ref{l:0y1} and~\ref{l:conexion}).

Take  $\Phi(x,y,\vec z)$ to be
\[\forall u :  \bigwedge_{i=1}^{k}  \Bigl( u\por U_{i}(x,y,\vec
0) = u \por U_{i}(x,y,\vec z) \Bigr) \ \longrightarrow \ u
\por x = u \por y.\]
This formula is factorable by Lemma~\ref{l:cota_inf}.
 Let $A, B\in \V$, and $a\in A$,
$b,d\in B$. First we prove that 
\[A\times B \models \Phi\bigl(\<a,b\>, \<a,d\>, [\vec 0, \vec
  1]\bigr)\]
Suppose that for some $\<u, v\>$ we have
\[A\times B \models \bigwedge_{i=1}^{k} \<u,v\> \por U_{i}(\<a,b\>,
\<a,d\>, [\vec 0, \vec 0]) = 
\<u,v\> \por U_{i}(\<a,b\>, \<a,d\>, [\vec 0, \vec 1]).\]
Then
\[ B \models\bigwedge_{i=1}^{k} v \por U_{i}(b,d,\vec 0) = v \por  U_{i}(b,d,\vec 1).\]
But the above equations in combination with~(\ref{eq:34}) produce
\[v \por b  = v \por d \] 
and hence
\[ \<u,v\> \por\<a,b\>= \<u,v\> \por  \<a,d\>.\]

Now suppose
\[A\times B \models \Phi\bigl(\<a,b\>, \<c,d\>, [\vec 0, \vec
  1]\bigr).\]
Since $\Phi$ is preserved by direct factors, we obtain $A \models
\Phi\bigl(a,c,\vec 0\bigr)$ and by inspection this is equivalent to
$\forall u : u \por a = u \por c$. 
%% Using this equation for
%% $u=a,c\por a,c$ we obtain:
%% \[a = a \por a  = a \por c = a \por c \por c  = c \por a \por c = c
%% \por a \por a = c \por a = c \por c = c,\]
Using this equation for
$u=a,c$ we obtain $a\men c$ and $c\men a$, therefore $a=c$.
\end{proof}

Using this new definition of factor congruences, we obtain:
\begin{theorem}
Let $\V$ be a semidegenerate variety of connected po-groupoids over
a finite language that satisfy~(\ref{eq:rel_semilatt}). Then the 
class of directly indecomposable algebras of $\V$ is axiomatizable by
a $\Pi_4$ sentence plus axioms for $\V$.
\end{theorem}

\bigskip

\begin{quote}
CIEM --- Facultad de Matem\'atica, Astronom\'{\i}a y F\'{\i}sica 
(Fa.M.A.F.) 

Universidad Nacional de C\'ordoba - Ciudad Universitaria

C\'ordoba 5000. Argentina.

\texttt{sterraf@famaf.unc.edu.ar}
\end{quote}

{\small
\begin{thebibliography}{10}
\bibitem{Ger1} \textsc{J. A. Gerhard,} \textit{The lattice of
  equational classes of idempotent semigroups,} J. Algebra \textbf{15}
  (1970): 195--224.
\bibitem{Ger2} \textsc{J. A. Gerhard,} \textit{Subdirectly irreducible
  idempotent semigroups,} Pacific J. Math. \textbf{39} (1971):
  669--676.
%% \bibitem{Hashimoto} \textsc{J. Hashimoto} \textit{On the direct
%%   product decomposition of partially ordered sets,} Ann. Math. \textbf{54} (1951):
%%   315--318.
\bibitem{5} \textsc{J. Kollar}, \textit{Congruences and one element
subalgebras,} Algebra univers. \textbf{9} (1979): 266--267.
\bibitem{Mace4}   \textsc{W. McCune,} \textit{Mace4 Reference Manual and Guide,} Tech. Memo ANL/MCS-TM-264, Mathematics and Computer Science Division, Argonne National Laboratory, Argonne, IL, August 2003. 
%% \bibitem{Otter}  \textsc{W. McCune,} \textit{Otter 3.3 Reference Manual,} Tech. Memo ANL/MCS-TM-263, Mathematics and Computer Science Division, Argonne National Laboratory, Argonne, IL, August 2003. 
\bibitem{p9m4} \textsc{W. McCune,} \textit{Prover9 and Mace4 Webpage,}  
  \texttt{http://www.cs.unm.edu/\~{ }mccune/prover9/}
\bibitem{4} \textsc{R. McKenzie, G. McNulty and W. Taylor},
\textit {Algebras, Lattices, Varieties,} Volume {\bf 1}, The Wadsworth
\& Brooks/Cole Math. Series, Monterey, California (1987).
\bibitem{DFC} \textsc{P. S\'anchez Terraf and D. Vaggione,} \textit{Varieties with Definable
Factor Congruences,} Trans. Amer. Math. Soc., to appear.
\bibitem{va0} \textsc{D. Vaggione}, \textit{$\V$ with factorable
congruences and $\V= \mathbf{I\Gamma}^a(\V_{DI})$ imply $\V$  is a
discriminator variety}.  Acta Sci. Math. 
\textbf{62} (1996): 359--368.
\bibitem{va5} \textsc{D.~Vaggione},
\textit{Varieties of shells}, Algebra univers. {\bf 36} (1996): 483--487.
\bibitem{indec-lat} \textsc{R. Willard}, \textit{A note on indecomposable
  lattices,} Algebra univers. {\bf 26} (1989): 257--258.
\bibitem{7} \textsc{R. Willard}, \textit{Varieties Having Boolean
Factor Congruences,} J. Algebra, \textbf{132} (1990): 130--153.
\end{thebibliography}
}
\end{document}